\documentclass[a4paper,12pt]{article}
\usepackage[latin1]{inputenc}
\usepackage{amsfonts}
\usepackage[dvips]{graphics}

\newcommand{\proxf}{\mathrm{prox}_f}
\newcommand{\proxH}{\mathrm{prox}_H}
\newcommand{\proxHs}{\mathrm{prox}_{H^\ast}}
\newtheorem{lem}{Lemma}
\newtheorem{theo}{Theorem}

\title{On a
generalization of the iterative soft-thresholding algorithm for the case of
non-separable penalty}
\author{Ignace Loris and Caroline Verhoeven\\
 Mathematics Department, Université Libre de Bruxelles,\\ CP
217,
Boulevard du Triomphe\\
B-1050 Bruxelles, Belgium\\
igloris@ulb.ac.be,cverhoev@ulb.ac.be}

\begin{document}
\maketitle

\begin{abstract}
An explicit algorithm for the minimization of an $\ell_1$
penalized least squares functional, with non-separable $\ell_1$
term, is proposed. Each step in the iterative algorithm
requires four matrix vector multiplications and a single simple
projection on a convex set (or equivalently thresholding).
Convergence is proven and a $1/N$ convergence rate is derived
for the functional. In the special case where the matrix in the
$\ell_1$ term is the identity (or orthogonal), the algorithm
reduces to the traditional iterative soft-thresholding
algorithm. In the special case where the matrix in the
quadratic term is the identity (or orthogonal), the algorithm
reduces to a gradient projection  algorithm for the dual
problem.

By replacing the projection with a simple proximity operator, other convex non-separable penalties than those based on an $\ell_1$-norm can be handled as well.
\end{abstract}

\section{Introduction}

Non-smooth minimization problems involving a sum of a quadratic
data misfit term and a non-smooth penalty term have received a
lot of attention in inverse problems and imaging in recent
years. In this note we are interested in finding the
minimizer $\hat x$ of the $\ell_1$ penalized least squares
functional $\mathcal{F}$:
\begin{equation}
\hat x=\arg\min_x\mathcal{F}(x) \qquad \mathrm{with}
\qquad\mathcal{F}(x)=\frac{1}{2}\|Kx-y\|^2+\lambda \|Ax\|_1,
\label{functional}
\end{equation}
by means of an iterative algorithm. Here $\|u\|^2=\sum_i u_i^2$
with $u_i\in\mathbb{R}$ and $\|w\|_1=\sum_i |w_i|$ ($w_i$ may
be an element of $\mathbb{R},\mathbb{R}^2,\ldots$ and $|w_i|$
stands for the Euclidean length of $w_i$; other choices of
$|w_i|$ are discussed in section \ref{discussionsection}). $K$
is a matrix mixing the variables in the quadratic data misfit
term and $A$ is a linear operator mixing the variables in the
penalty term. The quadratic term is convex and smooth, but the
penalty term $\|Ax\|_1$ is convex and non-smooth. We work in
a finite dimensional setting.

For the case where the non-smooth penalty term in (\ref{functional}) is simple
($A=1$) many algorithms have appeared in recent years. One of
the earliest (not necessarily the most efficient) is the
iterative soft-thresholding algorithm \cite{Daubechies2004b}
(see also section \ref{problemsection}). As $\ell_1$-norm
penalties promote sparsity, such algorithms are used in
`compressed sensing' \cite{Donoh2006} for finding a sparse
solution (up to noise level) of a large-scale under-determined
linear system. As problems in 2D and 3D imaging are large scale
problems, with many unknowns, such simple first-order iterative
algorithms can still be useful.

The principal difference of this paper with respect to
\cite{Daubechies2004b} is the presence of the matrix $A$ in the
penalty term. In image processing the total variation penalty,
which favors piece-wise constant images, is popular for its
ability to maintain sharp edges. The total variation penalty is
defined by the $\ell_1$-norm of the gradient of the unknown
($A=\mathrm{grad}$). It has mostly been studied for denoising
($K=1$) or for other special operators $K$ (e.g.
deconvolution).

Our aim here is to provide a simple iterative algorithm for the problem (\ref{functional}) with proven convergence (see
theorem \ref{theorem1}). We also desire an algorithm that is
fully explicit: each step in the proposed iteration only uses
four matrix-vector multiplications (one by $K,K^T,A$ and $A^T$)
and a simple projection on the $\ell_\infty$ ball (or
equivalently a single thresholding).

Although our main aim is to solve problem (\ref{functional}), we will formulate an algorithm and a convergence theorem for the more general problem:
\begin{equation}
\hat x=\arg\min_x\mathcal{F}(x) \qquad \mathrm{with}
\qquad\mathcal{F}(x)=\frac{1}{2}\|Kx-y\|^2+H(Ax),
\label{genfunctional}
\end{equation}
where $H$ is a convex function (we assume that the solution to (\ref{genfunctional}) exists).
For problem (\ref{genfunctional}) the projection operator mentioned before is replaced with the proximity operator $\proxHs$ of the convex conjugate $H^\ast$ and soft-thresholding is replaced with the proximity operator $\proxH$ of $H$. It is \emph{not} necessary to know the proximity operator of $H(A\,\cdot)$.

A second goal of the paper is to bridge the gap between the
well-known iterative soft-thresholding algorithm (used for the
special case $A=1$) and the general case $A\neq 1$ in problem (\ref{functional}). The
iterative soft-thresholding algorithm is well understood and
has a $1/N$ convergence rate for the decrease of the
functional. It is also the basis of an accelerated algorithm
with an improved $1/N^2$ rate of decrease of the functional
\cite{Beck.Teboulle2008,Nesterov1983a}. The averages of the first $N$ iterates of the
proposed generalized soft-thresholding algorithm are proven to have a $1/N$ rate on the
functional.

Our results differ from several existing algorithms for solving
(\ref{functional}) where each iteration step requires either
the solution of another (non-trivial) minimization problem, the
solution of a linear system, or a non-trivial projection on a
convex set. Our proposed algorithm may therefore be of use in
cases where the matrices involved ($K$ and $A$) have no special
structure that makes such sub-problems easily solvable (i.e.
not limited to deconvolution problems on regular grids, to
orthogonal matrices, etc.).

Iterative algorithms for the denoising case ($K=1$) can,
amongst others, be found in \cite{Chambolle2004,Chambolle2005}.
For general $K$, an algorithm that uses a smoothing parameter
is found in \cite{CHAN.GOLUB.ea1999}, an algorithm which needs
a projection on a non-trivial convex set is in \cite{eccv04}
and an algorithm which needs the solution of a non-trivial
sub-problem is in
\cite{Daubechies.Teschke.ea2007,Beck.Teboulle2009,Bredies2009}.
These are results for $A=\mathrm{grad}$ but this is not
essential in those algorithms.

Zhu and Chan \cite{Zhu.Chan2008} studied a primal-dual
formulation and a so-called `primal-dual hybrid gradient
descent' (PDHG) algorithm but concentrated on deconvolution.
Connections with (more general) algorithms for variational
inequalities were mentioned. This PDHG algorithm was placed in
a general framework for primal-dual algorithms in
\cite{Esser.Zhang.ea2010} and many interconnections can be
found there. The plethora of algorithms mentioned there still
require either the solution of a linear system (which may easy
in some special cases) or the minimization of a non-trivial
sub-problem. Applications to image recovery of an algorithm
that is an instance of the so-called alternating direction
method of multipliers, are tested in
\cite{Afonso.Bioucas-Dias.ea2010}.

Recently an \emph{explicit} algorithm was proposed in
\cite[equation 5.11]{Zhang.Burger.ea2011} with proven
convergence. No rate on the functional was given. That explicit
algorithm is different from the
one presented here. It does not reduce to the iterative
soft-thresholding algorithm when $A=1$. Another explicit
algorithm can also be derived using \cite[Eq.
74]{Chambolle.Pock2010} by the introduction of additional dual
variables.

It remains a subject of study what speed increase can be gained
(if any) from using an algorithm that solves a linear system at
every iteration. The derivation of an $\mathcal{O}(1/N^2)$
algorithm, if at all possible for this problem, would be more
interesting.  Our analysis and proof is inspired by
\cite{Popov1980,Pock.Cremers.ea2010} (who discuss a primal-dual
algorithm for another problem) and by \cite{Chambolle.Pock2010}. It is worth pointing out that no smoothing parameter is
introduced in the non-smooth part of the functional. The proposed
algorithm is not an iteratively reweighted least squares
algorithm.

\section{Mathematical tools}

We assume that the function $H(x)$ in (\ref{genfunctional}) and its convex conjugate $H^\ast(w)=\sup_x\langle w,x\rangle-H(x)$ are two proper, lower semi-continuous, convex functions on a finite dimensional real vector space and with image in $\mathbb{R}\cup \{+\infty\}$ \cite{Rockafellar1997}. For example in case of problem (\ref{functional}), $H(u)=\lambda \|u\|_1$ and therefore
\begin{equation}
H^\ast(w)=\left\{\begin{array}{lcl}
0&& \|w\|_\infty\leq \lambda\\
+\infty&& \|w\|_\infty>\lambda
\end{array}\right.\label{Hstarexample}
\end{equation}
such that
$\lambda\|u\|_1=H(u)=\max_w \langle w,u\rangle-H^\ast(w)=\sup_{\|w\|_\infty\leq \lambda}\langle w,u\rangle$.

The proximity operators \cite{Combettes.Pesquet2011} of the convex functions $H$ and $H^\ast$ are defined as:
\begin{equation}
\begin{array}{lcl}
\proxH(u)&=&\arg\min_x H(x)+\frac{1}{2}\|x-u\|^2\\
\proxHs(u)&=&\arg\min_w H^\ast(w)+\frac{1}{2}\|w-u\|^2.
\end{array}
\label{proxdef}
\end{equation}
Therefore, when $H(u)=\lambda \|u\|_1$ and $H^\ast$ is given by expression (\ref{Hstarexample}), we find that
$\proxHs(u)=P_\lambda(u)$, the projection on the $\ell_\infty$
ball of radius $\lambda$.
It has a simple explicit expression:
\begin{equation}
P_\lambda(u)=
\left\{
\begin{array}{ll}
\displaystyle\lambda\,\frac{u}{|u|} & |u|> \lambda\\
u & |u|\leq \lambda
\end{array}
\right.\label{proj}
\end{equation}
(applied component-wise).
On the other hand the proximity operator of $H(u)=\lambda\|u\|_1$ is the so-called soft-thresholding operator, $\proxH(u)=S_\lambda(u)$. It has the explicit expression:
\begin{equation}
S_\lambda(u)=
\left\{
\begin{array}{ll}
u-\frac{u}{|u|}\lambda &  |u|>\lambda\\
0 & |u|\leq \lambda\\
\end{array}
\right.\label{soft}
\end{equation}
(also applied component-wise).
Clearly, soft-thresholding $S_\lambda$ and projection $P_\lambda$ are connected by:
\begin{equation}
P_\lambda(u)=u-S_\lambda(u).\label{PSeq}
\end{equation}
In the formulas for $P_\lambda(u)$ and $S_\lambda(u)$ $u$ can be an element of $\mathbb{R},
\mathbb{R}^2,\ldots$ depending on context (in particular
$(Ax)_i\in\mathbb{R}^2$ when $A=\mathrm{grad}$ of a 2D image).
We shall use the same notation $S_\lambda, P_\lambda$ when
applied componentwise to a list of elements of $\mathbb{R},
\mathbb{R}^2,\ldots$

Proximity operators are Lipschitz-continuous mappings \cite{Combettes.Pesquet2011}:
\begin{equation}
\|\proxHs(u)-\proxHs(v)\|\leq\|u-v\|\qquad\qquad \forall u,v.\label{Lipschitzeq}
\end{equation}

The subdifferential $\partial H(x)=\{\gamma | H(y)\geq H(x)+\langle\gamma,y-x\rangle\quad \forall y\}$ of $H$ in $x$ can be characterized using the proximity operator of $H$. Indeed, from the definition (\ref{proxdef}) it follows that $u^+=\proxH(u^-)$ if and only if $0\in \partial H(u^+)+u^+-u^-$ or $u^--u^+\in \partial H(u^+)$. In other words, setting $u=u^--u^+$, we have that $u\in \partial H(u^+)$ if and only if $u^+=\proxH(u^++u)$.

Finally, it can also be shown that the proximity operator of $H$ and its dual $H^\ast$ are related by the following identity \cite{Moreau1965}:
\begin{equation}
\proxHs(u)+\proxH(u)=u,\label{proxsum}
\end{equation}
as already verified in equation (\ref{PSeq}) for the special case $\proxH=S_\lambda$ and $\proxHs=P_\lambda$.

We refer to \cite{Combettes.Pesquet2011} for a table of further examples and properties of proximity operators. The proximity operators $\proxHs=P_\lambda$ and $\proxH=S_\lambda$, that will be used for problem (\ref{functional}), i.e. for $H(u)=\lambda\|u\|_1$, have explicit expressions that are easy to implement.

\section{Variational equations and special cases}

\label{problemsection}

The variational equations of the minimization problem
(\ref{genfunctional}) are:
\begin{displaymath}
K^T(Kx-y)+A^T w=0,
\end{displaymath}
where $w$ is an element of the subdifferential of
$H(Ax)$. As mentioned before, this means that $Ax=\proxH(w+Ax)$ or equivalently, using
(\ref{proxsum}), that $w=\proxHs(w+Ax)$. The variational
equations corresponding to the problem (\ref{genfunctional}) are
therefore:
\begin{equation}
K^T (y-Kx)-A^T w=0\qquad\mathrm{and}\qquad
w=\proxHs\left(w+Ax\right).\label{vareq}
\end{equation}
The goal of this paper is to write an iterative algorithm that
converges to a solution of these equations. We assume that these equations have at least one solution $(\hat x,\hat w)$.

By using that
$H(Ax)=\sup_w\langle
Ax,w\rangle-H^\ast(w)$, the minimization problem (\ref{genfunctional}) can
also be written as a saddle-point problem
\begin{equation}
\min_x\max_w F(x,w),
\label{minmax}
\end{equation}
where we have set:
\begin{equation}
F(x,w)=\frac{1}{2}\|Kx-y\|^2+\langle Ax,w\rangle-H^\ast(w).
\label{Fdef}
\end{equation}
A saddle point $(\hat x,\hat w)$ of
(\ref{minmax}) is a point such that
\begin{equation}
F(\hat x,w)\leq F(\hat x,\hat w)\leq F(x,\hat w)\label{saddle}
\end{equation}
for all $x$ and $w$. For completeness,
we show in the next section that solutions $(\hat x,\hat w)$ of
equations (\ref{vareq}) are saddle-points of (\ref{minmax}). We
define the gap with respect to the saddle-point $(\hat x,\hat
w)$ by:
\begin{equation}
G(x,w)=F(x,\hat w)-F(\hat x,w).
\label{gapdef}
\end{equation}
It follows from (\ref{saddle}) that this gap is non-negative for
all $x$ and $w$.

In the special case $A=1$ the problem (\ref{genfunctional}) reduces to:
\begin{equation}
\min_x\frac{1}{2}\|Kx-y\|^2+H(x),
\end{equation}
for which a forward-backward splitting algorithm
\begin{equation}
x^{n+1}=\proxH\left(x^n+K^T(y-Kx^n)\right)\label{fbs}
\end{equation}
can be used. This algorithm converges for $\|K\|<\sqrt{2}$ \cite{ComWa2005}.
More specifically, the minimization problem with $A=1$ and $H(x)=\lambda\|x\|_1$:
\begin{equation}
\min_x\frac{1}{2}\|Kx-y\|^2+\lambda \|x\|_1
\end{equation}
can be solved by the iterative soft-thresholding
algorithm \cite{Daubechies2004b}:
\begin{equation}
x^{n+1}=S_\lambda\left(x^n+K^T(y-Kx^n)\right).\label{ista}
\end{equation}
Many other
algorithms exist. One feature of this algorithm is that,
as a consequence of the soft-thresholding, all the iterates
$x^n$ (not just the limit) have many exact zeros.

On the other hand, the problem
\begin{equation}
\min_x\frac{1}{2}\|x-g\|^2+\lambda \|Ax\|_1
\end{equation}
($K=1$, $y\rightarrow g$ in problem (\ref{functional})) can be solved by a gradient
projection algorithm:
\begin{equation}
w^{n+1}=P_\lambda\left(w^n+A(g-A^T w^n)\right)\label{projgrad}
\end{equation}
where $x^n=g-A^T w^n$, if $\|A\|<1$ (as is shown in \cite[eqn.
11]{Chambolle2005} for $A=\mathrm{grad}$). This is a special
case of the gradient projection algorithm that can be used for
minimization of a quadratic function over a convex set $C$:
$\min_{w\in C}\|g-A^Tw\|^2$. The quantities $Ax^n$ are not
sparse in every step, only in the limit will $Ax^n$ be sparse.

\section{Algorithm}

\label{algsection}

Writing the variational equations (\ref{vareq}) as
fixed-point equations:
\begin{equation}
\left\{
\begin{array}{lcl}
x&=&x+K^T (y-Kx)-A^T w\\
w&=&\proxHs\left(w+Ax\right),
\end{array}
\right.\label{fixeq}
\end{equation}
provides the usual ansatz for deriving iterative first order
algorithms for (\ref{genfunctional}). Here we choose to study the
iteration
\begin{equation}
\left\{
\begin{array}{lcl}
\bar x^{n+1}&=&x^n+K^T (y-Kx^n)-A^T w^n\\
w^{n+1}&=&\proxHs\left(w^n+A\bar x^{n+1}\right)\\
x^{n+1}&=&x^n+K^T (y-Kx^n)-A^T w^{n+1},
\end{array}
\right.   \label{alg}
\end{equation}
the fixed-point of which is a solution to the variational
equations (\ref{vareq}). Specifically, starting from
$(x^n,w^n)$ one does a gradient descent step on $F(x,w)$ in the
$x$-variable to arrive at $(\bar x^{n+1},w^n)$, followed by a
proximal ascent step in the $w$ variable to compute
$w^{n+1}$. Finally one does a gradient descent step in
$(x^n,w^{n+1})$ to arrive at $(x^{n+1},w^{n+1})$. This
algorithm can therefore be interpreted as a `predict-correct'
algorithm for the saddle-point problem (\ref{minmax}). On the
other hand the algorithm (\ref{alg}) can equivalently be
written in a `pseudo-implicit' form as:
\begin{equation}
\left\{
\begin{array}{lcl}
\bar x^{n+1}&=&x^{n+1}-A^T (w^n-w^{n+1})\\
w^{n+1}&=&\proxHs\left(w^n+A\bar x^{n+1}\right)\\
x^{n+1}&=&x^n+K^T (y-Kx^n)-A^T w^{n+1}.
\end{array}
\right.
\end{equation}
This form is useful for proving convergence.

Writing the algorithm (\ref{alg}) as:
\begin{equation}
\left\{
\begin{array}{lcl}
g^{n+1}&=&x^n+K^T (y-Kx^n)\\
w^{n+1}&=&\proxHs\left(w^n+A (g^{n+1}-A^T w^n)\right)\\
x^{n+1}&=&g^{n+1}-A^T w^{n+1},
\end{array}
\right.\label{algalt}
\end{equation}
leads to the interpretation of a gradient descent step on the
quadratic part of the functional, followed by a single step in
a dual variable (compare with (\ref{projgrad})) starting from
the previous dual variable $w^n$.

In the next section we show that the proposed algorithm
(\ref{alg}) converges to a solution of the fixed-point
equations (\ref{fixeq}), i.e. to a saddle-point of the min-max
problem (\ref{minmax}) and to a minimizer of the functional
(\ref{genfunctional}). Under some additional condition on $H$ we also derive a convergence rate estimate for
the functional $\mathcal{F}$ in the average of the iterates.

For the special case when $AA^T=A^TA=1$, the second
line of algorithm (\ref{alg}) reduces to:
\begin{displaymath}
w^{n+1}=\proxHs\left(A(x^n+K^T (y-Kx^n))\right)
\end{displaymath}
which implies:
\begin{displaymath}
x^{n+1}=x^n+K^T (y-Kx^n)-A^T \proxHs\left(A(x^n+K^T
(y-Kx^n))\right).
\end{displaymath}
Using $\proxH(u)=u-\proxHs(u)$, one has:
\begin{displaymath}
x^{n+1}= A^T \proxH(\left(A(x^n+K^T (y-K x^n))\right).
\end{displaymath}
This is the forward-backward splitting algorithm (\ref{fbs}) for the
variable $Ax$ and the operator $KA^T$. In particular, for  $H=\lambda\|\cdot\|_1$, the algorithm (\ref{alg}) reduces to the iterative soft-thresholding algorithm (\ref{ista}) when $AA^T=A^TA=1$. Similarly, when $K$ is
orthogonal and $H(\cdot)=\lambda\|\cdot\|_1$, then the algorithm (\ref{alg}) reduces to
\begin{displaymath}
\left\{
\begin{array}{lcl}
w^{n+1}&=&P_\lambda\left(w^n+A(K^Ty-A^T w^n)\right)\\
x^{n+1}&=&K^T y-A^T w^{n+1}
\end{array}
\right.
\end{displaymath}
which is the gradient projection algorithm (\ref{projgrad}) for
the data $g=K^Ty$.

\section{Convergence}
\label{proofsection}

We will prove convergence of algorithm (\ref{alg}).
\begin{lem} If $w^+=\proxHs(w^-+\Delta)$ then
\begin{equation}
\|w-w^+\|^2\leq \|w-w^-\|^2-\|w^--w^+\|^2-2\langle w-w^+,\Delta\rangle+2H^\ast(w)-2H^\ast(w^+)\label{lemmaproxeq}
\end{equation}
for all $w$. \label{lemmaprox}
\end{lem}
Proof: If $w^+=\proxHs(w^-+\Delta)=\arg\min_w H^\ast(w)+\frac{1}{2}\|w-(w^-+\Delta)\|^2$ then $w^-+\Delta-w^+\in\partial H^\ast(w^+)$. We then have for all $w$ that $H^\ast(w) \geq H^\ast(w^+)+\langle \gamma,w-w^+\rangle $ if $\gamma\in\partial H^\ast(w^+)$ and therefore:
\begin{displaymath}
\begin{array}{lcl}
H^\ast(w) &\geq& H^\ast(w^+)+\langle w^-+\Delta-w^+,w-w^+\rangle\\
 & = & H^\ast(w^+)+\langle \Delta,w-w^+\rangle+\langle w^--w^+,w-w^+ \rangle\\
 & = & H^\ast(w^+)+\langle \Delta,w-w^+\rangle+\frac{1}{2}\|w^--w^+\|^2+\frac{1}{2}\|w-w^+\|^2-\frac{1}{2}\|w^--w\|^2
\end{array}
\end{displaymath}
which gives (\ref{lemmaproxeq}).\hfill$\Box$

\begin{lem} If $x^+=x^-+\Delta$, then
\begin{equation}
\|x-x^+\|^2=\|x-x^-\|^2-\|x^--x^+\|^2-2\langle x-x^+,\Delta\rangle
\end{equation} for all $x$.
\label{lemmax}
\end{lem}
For completeness we show that a solution of the variational
equations is a saddle-point of (\ref{minmax}). This implies
that the gap $G(x,w)$ with respect to the fixed-point $(\hat x,
\hat w)$ is always non-negative.
\begin{lem}
If $(\hat x,\hat w)$ satisfies the fixed-point equations
(\ref{fixeq}), then
\begin{equation}
F(\hat x,w)\leq F(\hat x,\hat w)\leq F(x,\hat w)
\end{equation}
and hence
\begin{equation}
G(x,w)\equiv F(x,\hat w)-F(\hat x,w)\geq 0
\end{equation}
for all $x$ and $w$. \label{lemmagap}
\end{lem}
Proof: The first inequality $F(\hat x,w)\leq F(\hat x,\hat w)$
comes down to showing that $0\leq\langle A\hat x,\hat
w-w\rangle+H^\ast(w)-H^\ast(\hat w)$ for all $w$. This follows
immediately from choosing $w^+=w^-=\hat w$ and $\Delta=A\hat
x$ in lemma \ref{lemmaprox}.

The second inequality $F(\hat x,\hat w)\leq F(x,\hat w)$ can be
written as:
\begin{displaymath}
0\leq\frac{1}{2}\|Kx-y\|^2-\frac{1}{2}\|K\hat x-y\|^2+\langle A
(x-\hat x),\hat w\rangle\qquad \forall x.
\end{displaymath}
To show this we choose $x^+=x^-=\hat x$ and $\Delta=K^T(y-K\hat
x)-A^T\hat w$ in lemma \ref{lemmax} to find:
\begin{displaymath}
\begin{array}{lcl}
0&=&-2\langle x-\hat x,K^T(y-K\hat x)-A^T\hat w\rangle\\
&=&  -2\langle K(x-\hat x),y-K\hat x\rangle+2\langle x-\hat x,A^T\hat w\rangle\\
&=& -\|K(x-\hat x)\|^2-\|y-K\hat x\|^2+\|Kx-y\|^2+2\langle x-\hat x,A^T\hat w\rangle\\
\end{array}
\end{displaymath}
for all $x$, or
\begin{displaymath}
\|K(x-\hat x)\|^2= \|Kx-y\|^2-\|K\hat x-y\|^2+2\langle A
(x-\hat x),\hat w\rangle,
\end{displaymath}
which is a slightly stronger result than needed.\hfill$\Box$

The gap $G(x,w)$ equals:
\begin{equation}
G(x,w)=\frac{1}{2}\|K(\hat x-x)\|^2+\langle \hat w-w,A\hat x\rangle+H^\ast(w)-H^\ast(\hat w) \label{temp2}
\end{equation}
as can be verified from its definition (and lemma
\ref{lemmax}). The sum of the last three terms on the right hand side is non-negative, so
\begin{equation}
G(x,w)\geq \frac{1}{2}\|K(\hat x-x)\|^2. \label{temp5}
\end{equation}
The gap $G(x,w)$ is not a
measure of closeness of $(x,w)$ to a saddle-point $(\hat x,\hat
w)$ as  $G(x,w)=0$ does not imply that $(x,w)$ is a saddle
point.

\begin{lem} If $(x^n,w^n)$ are given by iteration (\ref{alg}) then
\begin{displaymath}
\begin{array}{lcl}
\|x-x^{n+1}\|^2+\|w-w^{n+1}\|^2&\leq&\|x-x^n\|^2+\|w-w^n\|^2\\
&& -\|x^n-x^{n+1}\|^2-\|w^n-w^{n+1}\|^2\\
& &-\|K(x-x^n)\|^2+\|K(x^n-x^{n+1})\|^2\\
& & -\|A^T(w-w^n)\|^2\\
& & +\|A^T(w^n-w^{n+1})\|^2+\|A^T(w-w^{n+1})\|^2  \\
 & & -2 \left(F(x^{n+1},w)-F(x,w^{n+1})\right)
\end{array}
\end{displaymath}
for all $x$ and $w$.\label{lemma5}
\end{lem}
Proof: From lemmas \ref{lemmaprox} and \ref{lemmax} we find:
\begin{displaymath}
\begin{array}{lcl}
\|w-w^{n+1}\|^2&\leq& \|w-w^n\|^2-\|w^n-w^{n+1}\|^2-2\langle w-w^{n+1},A\bar x^{n+1}\rangle\\
&&\qquad\qquad\qquad\qquad\qquad\qquad\qquad+2H^\ast(w)-2H^\ast(w^{n+1})\\
\|x-x^{n+1}\|^2&=&\|x-x^n\|^2-\|x^n-x^{n+1}\|^2-2\langle x-x^{n+1},K^T(y-Kx^n)-A^Tw^{n+1}\rangle
\end{array}
\end{displaymath}
which together yield:
\begin{displaymath}
\begin{array}{lcl}
\|x-x^{n+1}\|^2+\|w-w^{n+1}\|^2&\leq& \|x-x^n\|^2-\|x^n-x^{n+1}\|^2\\
&&+\|w-w^n\|^2-\|w^n-w^{n+1}\|^2\\
&& -2\langle w-w^{n+1},A\bar x^{n+1}\rangle\\
&&-2\langle x-x^{n+1},K^T(y-Kx^n)-A^Tw^{n+1}\rangle\\
&& +2H^\ast(w)-2H^\ast(w^{n+1})
\end{array}
\end{displaymath}
As (\ref{alg}) implies $\bar x^{n+1}=x^{n+1}-A^T(w^n-w^{n+1})$,
this can be written as:
\begin{displaymath}
\begin{array}{lcl}
\|x-x^{n+1}\|^2+\|w-w^{n+1}\|^2&\leq&\|x-x^n\|^2-\|x^n-x^{n+1}\|^2\\
&& +\|w-w^n\|^2-\|w^n-w^{n+1}\|^2\\
&&-2\langle w-w^{n+1},A(x^{n+1}-A^T(w^n-w^{n+1}))\rangle\\
&&-2\langle x-x^{n+1},K^T(y-Kx^n)-A^Tw^{n+1}\rangle\\
&& +2H^\ast(w)-2H^\ast(w^{n+1})
\end{array}.
\end{displaymath}
The two $\langle w^{n+1},Ax^{n+1}\rangle$ terms cancel:
\begin{displaymath}
\begin{array}{lcl}
\|x-x^{n+1}\|^2+\|w-w^{n+1}\|^2&\leq&\|x-x^n\|^2-\|x^n-x^{n+1}\|^2\\
&&+\|w-w^n\|^2-\|w^n-w^{n+1}\|^2\\
&& +2\langle A^T(w-w^{n+1}),A^T(w^n-w^{n+1})\rangle\\
&& -2\langle K(x-x^{n+1}),y-Kx^n\rangle\\
&&-2\langle Ax^{n+1},w\rangle+2\langle x,A^Tw^{n+1}\rangle\\
&& +2H^\ast(w)-2H^\ast(w^{n+1}).
\end{array}
\end{displaymath}
Now, by using the equalities:
\begin{displaymath}
\begin{array}{lcl}
2\langle A^T(w-w^{n+1}),A^T(w^n-w^{n+1})\rangle &=&-\|A^T(w-w^n)\|^2+\|A^T(w-w^{n+1})\|^2\\
&&\qquad\qquad\qquad\qquad\qquad+\|A^T(w^n-w^{n+1})\|^2\\
 -2\langle K(x-x^{n+1}),y-Kx^n\rangle &=&\|Kx-y\|^2-\|Kx^{n+1}-y\|^2-\|K(x-x^n)\|^2\\
 &&\qquad\qquad\qquad\qquad\qquad\qquad +\|K(x^n-x^{n+1})\|^2\\
-2\langle Ax^{n+1},w\rangle+2\langle x,A^Tw^{n+1}\rangle &=& 2F(x,w^{n+1})-2F(x^{n+1},w)-\|Kx-y\|^2\\
&&\qquad+\|Kx^{n+1}-y\|^2+2H^\ast(w^{n+1})-2H^\ast(w),
\end{array}
\end{displaymath}
the previous inequality reduces to:
\begin{displaymath}
\begin{array}{lcl}
\|x-x^{n+1}\|^2+\|w-w^{n+1}\|^2&\leq&\|x-x^n\|^2-\|x^n-x^{n+1}\|^2\\
&& +\|w-w^n\|^2-\|w^n-w^{n+1}\|^2\\
&& -\|A^T(w-w^n)\|^2+\|A^T(w^n-w^{n+1})\|^2\\
&&+\|A^T(w-w^{n+1})\|^2\\
&& -\|K(x-x^n)\|^2+\|K(x^n-x^{n+1})\|^2\\
&&+2F(x,w^{n+1})-2F(x^{n+1},w),
\end{array}
\end{displaymath}
which is the desired result.\hfill$\Box$

\begin{theo}
Let $\|K\|<\sqrt{2}$ and $\|A\|<1$. If the equations (\ref{vareq}) have a solution and the sequence
$(x^n,w^n)$ is defined by the iteration
\begin{equation}
\left\{
\begin{array}{lcl}
\bar x^{n+1}&=&x^n+K^T (y-Kx^n)-A^T w^n\\
w^{n+1}&=&\proxHs\left(w^n+A\bar x^{n+1}\right)\\
x^{n+1}&=&x^n+K^T (y-Kx^n)-A^T w^{n+1},
\end{array}
\right.   \label{theoalg}
\end{equation}
then:
\begin{enumerate}
\item the sequence $(x^n,w^n)$ converges to a solution
    $(x^\dagger,w^\dagger)$ of the variational equations
    (\ref{vareq}) thereby providing a minimizer of
    (\ref{genfunctional}) and a saddle point of
    (\ref{minmax}),

\item the average of the first $N$ iterates $(\tilde
    x^N,\tilde w^N)=\sum_{i=1}^N(x^i,w^i)/N$, converges to
    the saddle-point $(x^\dagger,w^\dagger)$ and there exists a constant $C_1\geq0$ independent of $N$ such that:
\begin{equation}
F(\tilde x^{N},w)-F(x,\tilde w^N)\leq \frac{\|x-x^0\|^2+\|w-w^0\|^2+C_1}{2N}
\label{th1}
\end{equation}
for all $x,w$, (with $C_1=0$ if $\|K\|\leq 1$), in particular:
\begin{equation}
0\leq G(\tilde x^N,\tilde w^N)\leq \frac{\|x^\dagger-x^0\|^2+\|w^\dagger-w^0\|^2+C_1}{2N}.
\label{th2}
\end{equation}
\item \label{point3}If the dual variable $w$ is bounded ($H^\ast(w)=+\infty$ for $\|w\|>R$ for some $R>0$), there exists a constant $C_2$ independent of $N$ such that:
\begin{equation}
0\leq \mathcal{F}(\tilde x^N)-\mathcal{F}(x^\dagger)\leq C_2/N\qquad\qquad\forall N.
\label{th3}
\end{equation}
\end{enumerate}
\label{theorem1}
\end{theo}
Proof: i) Let $(\hat x, \hat w)$ be a saddle point of
(\ref{minmax}). From lemma \ref{lemma5} we find:
\begin{displaymath}
\begin{array}{lcl}
\|\hat x-x^{n+1}\|^2+\|\hat w-w^{n+1}\|^2&\leq&\|\hat x-x^n\|^2-\|K(\hat x-x^n)\|^2\\
&&+\|\hat w-w^n\|^2-\|A^T(\hat w-w^n)\|^2\\
&&-\|x^n-x^{n+1}\|^2+\|K(x^n-x^{n+1})\|^2\\
& & -\|w^n-w^{n+1}\|^2+\|A^T(w^n-w^{n+1})\|^2 \\
& &+\|A^T(\hat w-w^{n+1})\|^2-\|K(\hat x-x^{n+1})\|^2 \\
\end{array}
\end{displaymath}
where we have used relation (\ref{temp5}) to set $2F(\hat
x,w^{n+1})-2F(x^{n+1},\hat w)=-2G(x^{n+1},w^{n+1})\leq
-\|K(\hat x-x^{n+1})\|^2$. Using the inequality
\begin{displaymath}
\begin{array}{lcl}
-\|K(\hat x-x^n)\|^2-\|K(\hat
x-x^{n+1})\|^2&=&-\frac{1}{2}\|K(\hat x-x^n)+K(\hat
x-x^{n+1})\|^2\\
&& \qquad\qquad-\frac{1}{2}\|K(\hat x-x^n)-K(\hat x-x^{n+1})\|^2\\
 &\leq & -\frac{1}{2}\|K(x^n-x^{n+1})\|^2\rule{0mm}{5mm},
\end{array}
\end{displaymath}
we find:
\begin{displaymath}
\begin{array}{lcl}
\|\hat x-x^{n+1}\|^2+\|\hat w-w^{n+1}\|^2&\leq&\|\hat x-x^n\|^2+\|\hat w-w^n\|^2-\|A^T(\hat w-w^n)\|^2\\
&&-\|x^n-x^{n+1}\|^2+\frac{1}{2}\|K(x^n-x^{n+1})\|^2\\
& & -\|w^n-w^{n+1}\|^2+\|A^T(w^n-w^{n+1})\|^2 \\
& &+\|A^T(\hat w-w^{n+1})\|^2.
\end{array}
\end{displaymath}

As we assume that $\|K\|<\sqrt{2}$ and $\|A\|<1$ we can
introduce regular square matrices $L$ and $B$ by
$L^TL=1-\frac{1}{2}K^TK$ and $B^TB=1-AA^T$ and deduce:
\begin{displaymath}
\begin{array}{lcl}
\|\hat x-x^{n+1}\|^2+\|B(\hat w-w^{n+1})\|^2&\leq&\|\hat x-x^n\|^2+\|B(\hat w-w^n)\|^2\\
&& -\|L(x^n-x^{n+1})\|^2-\|B(w^n-w^{n+1})\|^2.
\end{array}
\end{displaymath}
Summing from $N$ to $M\geq N$, one also finds:
\begin{equation}
\begin{array}{lcl}
\|\hat x-x^{M+1}\|^2+\|B(\hat w-w^{M+1})\|^2&\leq&\|\hat x-x^N\|^2+\|B(\hat w-w^N)\|^2\\
&& \!\!\!\!-\sum_{n=N}^M\|L(x^n-x^{n+1})\|^2+\|B(w^n-w^{n+1})\|^2.\label{temp1}
\end{array}
\end{equation}

As $B$ is invertible, it follows that the sequence $(x^n,w^n)$
is bounded. Hence there is a convergent subsequence
$(x^{n_j},w^{n_j})\stackrel{j\rightarrow\infty}{\rightarrow}(x^\dagger,w^\dagger)$
(the same subsequence for $x^n$ and $w^n$). It also follows
from inequality (\ref{temp1}) that:
\begin{equation}
\sum_{n=N}^M\|L(x^n-x^{n+1})\|^2+\|B(w^n-w^{n+1})\|^2\leq\|\hat
x-x^N\|^2+\|B(\hat w-w^N)\|^2.\label{temp4}
\end{equation}
Hence $\|L(x^n-x^{n+1})\|^2$ and $\|B(w^n-w^{n+1})\|^2$ tend to
zero for large $n$, which implies that $\|x^n-x^{n+1}\|$ and
$\|w^n-w^{n+1}\|$ tend to zero. It follows that the subsequence
$(x^{n_j+1},w^{n_j+1})$ also converges to
$(x^\dagger,w^\dagger)$ and, by continuity of $\proxHs$, that $(x^\dagger,w^\dagger)$
satisfies the fixed-point equations (\ref{fixeq}). We can
therefore choose $(\hat x, \hat w)=(x^\dagger,w^\dagger)$ in
relation (\ref{temp1}) to find:
\begin{equation}\label{temp6}
\|x^\dagger-x^{M+1}\|^2+\|B(w^\dagger-w^{M+1})\|^2\leq\|x^\dagger-x^N\|^2+\|B(w^\dagger-w^N)\|^2
\end{equation}
for all $M\geq N$. As there is a convergent subsequence of
$(x^n,w^n)$, the right hand side of this expression can be made
arbitrarily small for large enough $N$ ($N=n_j$ for some $j$).
Hence the left hand side will be arbitrarily small for all $M$
larger than this $N$. This proves convergence of the whole
sequence $(x^n,w^n)$ to $(x^\dagger,w^\dagger)$.

ii) As $(x^n,w^n)\stackrel{n\rightarrow\infty}{\rightarrow}
(x^\dagger,w^\dagger)$, the Césaro averages $(\tilde
    x^N,\tilde w^N)=\sum_{n=1}^N(x^n,w^n)/N$ also converge
to $(x^\dagger,w^\dagger)$. It follows from lemma \ref{lemma5}
that:
\begin{equation}
\begin{array}{lcl}
2 \left(F(x^{n+1},w)-F(x,w^{n+1})\right)&\leq&\|x-x^n\|^2+\|B(w-w^n)\|^2\\
  &&-\|x-x^{n+1}\|^2-\|B(w-w^{n+1})\|^2\\
  && -\|x^n-x^{n+1}\|^2+\|K(x^n-x^{n+1})\|^2.
\end{array}\label{temp3}
\end{equation}
Then, using convexity, one finds:
\begin{displaymath}
\begin{array}{lcl}
F(\tilde x^N,w)-F(x,\tilde w^N)&\leq& \displaystyle\frac{1}{N}\sum_{n=0}^{N-1} F(x^{n+1},w)-F(x,w^{n+1})\\
& \stackrel{(\ref{temp3})}{\leq} & \displaystyle\frac{1}{2N}\sum_{n=0}^{N-1} \|x-x^n\|^2+\|B(w-w^n)\|^2\\
&& \displaystyle\qquad\qquad-\|x-x^{n+1}\|^2-\|B(w-w^{n+1})\|^2\\
&& \displaystyle \qquad\qquad\qquad -\|x^n-x^{n+1}\|^2+\|K(x^n-x^{n+1})\|^2\\
&=&\displaystyle\frac{1}{2N}\left(\rule{0mm}{8mm}\|x-x^0\|^2+\|B(w-w^0)\|^2\right.\\
&& \displaystyle\qquad\qquad-\|x-x^{N}\|^2-\|B(w-w^{N})\|^2\\
&& \displaystyle \qquad\qquad\qquad \left.-\sum_{n=0}^N\|x^n-x^{n+1}\|^2-\|K(x^n-x^{n+1})\|^2\right)\\
& \leq & \displaystyle\frac{1}{2N}\left(\rule{0mm}{7mm} \|x-x^0\|^2+\|B(w-w^0)\|^2\right.\\
&& \qquad\qquad \displaystyle\left.-\sum_{n=0}^N\|x^n-x^{n+1}\|^2-\|K(x^n-x^{n+1})\|^2\right).
\end{array}
\end{displaymath}
If $\|K\|\leq 1$ the summation on the right hand side can be
dropped outright. If $1<\|K\|<\sqrt{2}$, the
$\|x^n-x^{n+1}\|^2$ terms can be dropped, and the sum
$\sum_{n=0}^N\|K(x^n-x^{n+1})\|^2$ can be bounded by the series
$\sum_{n=0}^\infty\|K(x^n-x^{n+1})\|^2$. The latter converges
as a consequence of relation (\ref{temp4}) and the regularity
of the matrix $L$. Therefore, a constant $C_1$ independent of
$N$ (and equal to $0$ in case $\|K\|\leq 1$), can be introduced
such that
\begin{displaymath}
F(\tilde x^N,w)-F(x,\tilde w^N) \leq
\displaystyle\frac{\|x-x^0\|^2+\|w-w^0\|^2+C_1}{2N}
\end{displaymath}
(where we have used $\|B\|\leq 1$).

Relation (\ref{th2}) follows from choosing
$(x,w)=(x^\dagger,w^\dagger)$ in equation (\ref{th1}).

iii) In case $H^\ast$ is such that the dual variable $w$ is bounded ($H^\ast=+\infty$ for all $w$ with $\|w\|>R$, for some $R>0$), we find
\begin{displaymath}
\begin{array}{lcl}
0\leq \mathcal{F}(\tilde x^N)-\mathcal{F}(x^\dagger)&=& \displaystyle\mathcal{F}(\tilde x^N)-F(x^\dagger,w^\dagger)\\
&\stackrel{\mathrm{lemma}\ \ref{lemmagap}}{\leq}&\displaystyle\mathcal{F}(\tilde x^N)-F(x^\dagger,\tilde w^N)     \\
&=&\displaystyle\sup_w F(\tilde x^N,w)-F(x^\dagger,\tilde w^N)     \\
&=&\displaystyle\max_{\|w\|\leq R} F(\tilde x^N,w)-F(x^\dagger,\tilde w^N) \\
 & \stackrel{(\ref{th1})}{\leq} & \displaystyle\max_{\|w\|\leq R}\frac{\|x^\dagger-x^0\|^2+\|w-w^0\|^2+C_1}{2N} \\
 &=& \displaystyle C_2/N
\end{array}
\end{displaymath}
which proves relation (\ref{th3}).
\hfill$\Box$

\section{Discussion}
\label{discussionsection}

In case the conditions $\|K\|<\sqrt{2}$ and/or $\|A\|<1$ are not satisfied, it is possible to rescale the matrices, the data $y$ and the variable $w$ to write a convergent algorithm. In the special case of functional (\ref{functional}) it suffices to rewrite the problem equivalently as:
\begin{displaymath}
\min_x \frac{1}{2}\|(\sqrt{\tau} K) x- (\sqrt{\tau} y)\|^2+\frac{\tau\lambda}{\sqrt{\sigma}} \|(\sqrt{\sigma}A)x\|_1
\end{displaymath}
and use algorithm (\ref{alg}) with $\proxHs=P_{\tau\lambda/\sqrt{\sigma}}$, for the matrices $\sqrt{\tau}K$, $\sqrt{\sigma}A$ and the data $\sqrt{\tau} y$.
Renaming $w\leftarrow \sqrt{\sigma}w/\tau$ and using the scaling property $P_{\tau\lambda/\sqrt{\sigma}}(u)=\tau/\sqrt{\sigma} P_\lambda(\sqrt{\sigma}u/\tau)$ one finds the
following iteration:
\begin{equation}
\left\{
\begin{array}{lcl}
\bar x^{n+1}&=&x^n+\tau K^T (y-Kx^n)-\tau A^T w^n\\
w^{n+1}&=&P_\lambda\left(w^n+\sigma/\tau A\bar x^{n+1}\right)\\
x^{n+1}&=&x^n+\tau K^T (y-Kx^n)-\tau A^T w^{n+1}
\end{array}
\right.   \label{scaledalg}
\end{equation}
for problem (\ref{functional}). Now step size parameters $\sigma,\tau>0$ should satisfy
$\tau<2/\|K^TK\|$ and $\sigma<1/\|AA^T\|$. In this case the bound (\ref{th3}) is valid as the dual variable $w$ is bounded (an element of the $\ell_\infty$ ball of radius $\lambda$).

For the general case (\ref{genfunctional}), the scaled version of the algorithm can be derived in a similar fashion. It takes the form
\begin{equation}
\left\{
\begin{array}{lcl}
\bar x^{n+1}&=&x^n+\tau K^T (y-Kx^n)-\tau A^T w^n\\
w^{n+1}&=&\mathrm{prox}_{\frac{\sigma}{\tau}H^\ast}\left(w^n+\sigma/\tau A\bar x^{n+1}\right)\\
x^{n+1}&=&x^n+\tau K^T (y-Kx^n)-\tau A^T w^{n+1}
\end{array}
\right.   \label{scaledalg2}
\end{equation}
with $\tau<2/\|K^TK\|$ and $\sigma<1/\|AA^T\|$.
Here we have used that $(\tau H(\cdot/\sqrt{\sigma}))^\ast=\tau H^\ast(\sqrt{\sigma}/\tau\,\cdot)$ and $\proxf(\alpha u)=\alpha\, \mathrm{prox}_{\alpha^{-2}f(\alpha\cdot)}(u)$ for $\alpha>0$.

It can be verified numerically that the functional
$\mathcal{F}(x^n)$ does not necessarily decrease monotonically as a
function of $n$ (this can be shown to hold in the special case
$A=1$ and $H(u)=\lambda\|u\|_1$, see \cite{Beck.Teboulle2008}). The gap function $G(x^n,w^n)$
does not decrease monotonically as a function of $n$ either.
The error between $(x^n,w^n)$ and $(x^\dagger,w^\dagger)$
decreases monotonically as a function of $n$ in the norm
$(\|x\|^2+\|Bw\|^2)^\frac{1}{2}$. This is a consequence of relation (\ref{temp6}).

The condition $\|A\|<1$ used in the proof of convergence
excludes the case $A=1$. Nevertheless, the proof of convergence
in theorem \ref{theorem1} can be slightly adapted to cover the
case $A=1$ as well.

The strength of algorithm (\ref{alg}) lies in the fact that only $\proxHs$ is needed and not $\mathrm{prox}_{H(A\cdot)}$. $\proxHs$ may have a simple expression whereas the proximity operator of $H(A\,\cdot)$ may not. In particular, for $H=\lambda\|\cdot\|_1$, one has expressions (\ref{proj}) and (\ref{soft}) for $\proxHs$ and $\proxH$. One can also find closed-form expressions for $\proxHs$ when  $|\cdot|$ (used in the expression $\|Ax\|_1=\sum_i|(Ax)_i|$) refers to the $1$- or $\infty$-norms instead of the 2-norm.

We believe the proposed algorithm (\ref{scaledalg}), its connection with the
traditional iterative soft-thresholding algorithm and its proof
of convergence are new. The combination of a gradient step with
the dual algorithm (\ref{projgrad}) has been proposed several
times already \cite{Beck.Teboulle2009,Bredies2009}; as such
that would not be an explicit algorithm as it requires
infinitely many dual iterations in each outer iteration. Here
we have shown convergence in the case when just one dual step
is made in each iteration. The series of algorithms discussed
in \cite{Esser.Zhang.ea2010,Zhang.Burger.ea2011} mostly make
use of a non-explicit step in the iteration, or of the solution
of a linear system at every iteration. These existing
algorithms are often special cases of more general methods. The
explicit algorithm in \cite{Zhang.Burger.ea2011} is also
different. In \cite[Eq. 74]{Chambolle.Pock2010} the authors
propose another explicit method using additional dual
variables.

In \cite{Korpelevich1976,Korpelevich1977} Korpelevich introduced an extragradient algorithm for the solution of a general saddle point problem. It relies on updating two copies of primal and dual variables (say $(x^n,w^n)$ and $(\bar x^n,\bar w^n)$), and combining the one with the gradient in the other point. No distinction is made there between primal and dual variables, as we do here.

In \cite{Nemirovski2005} a proximal point algorithm is introduced that reduces to the Korpelevich algorithm in a special case. Moreover it was shown that the iterates of that algorithm converge `ergodically' with rate $1/N$ on the primal objective function (i.e. the Césaro means decrease the functional with rate $1/N$). It is however assumed there that both primal and dual variables are bounded (saddle point problem on a compact domain). This is a major difference with our result, where we only require that the dual variable $w$ be bounded in order to derive a $1/N$ bound on the functional (point \ref{point3} of theorem \ref{theorem1}); the primal variable $x$ is unbounded.

We did not try to extend the convergence proof to an infinite
dimensional setting, as was done in \cite{Daubechies2004b} for
the iterative soft-thresholding algorithm. The most useful
example of problem (\ref{functional}) is perhaps the case where
$A=\mathrm{grad}$ (total variation penalty), but this operator
is unbounded in the infinite dimensional case.

Algorithm (\ref{scaledalg}) for problem (\ref{functional}) can also be used for signal recovery under  analysis style sparsity requirements \cite{Nama2011}: finding $x$ with many $(Ax)_i$ equal to zero for a given frame operator $A$.
Another application of algorithm (\ref{scaledalg}) is solving a linear inverse problem while imposing group sparsity (possibly with overlapping groups) \cite{Yuan2006}. In this case the matrix $A$ is chosen in such a way that $\|Ax\|_1=\sum_{k\in\mathrm{groups}} |(x_i)_{i\in\mathrm{group}_k}|$, i.e $A$ has a single $1$ on each row (and all other elements are zero). Columns may have more than a single nonzero entry (this would correspond to overlapping groups). In this expression $|\cdot|$ is again the Euclidean norm of a vector.

\section{Numerical example of total variation minimization}
\label{examplesection}

\begin{figure}
\centering\resizebox{\textwidth}{!}{\includegraphics{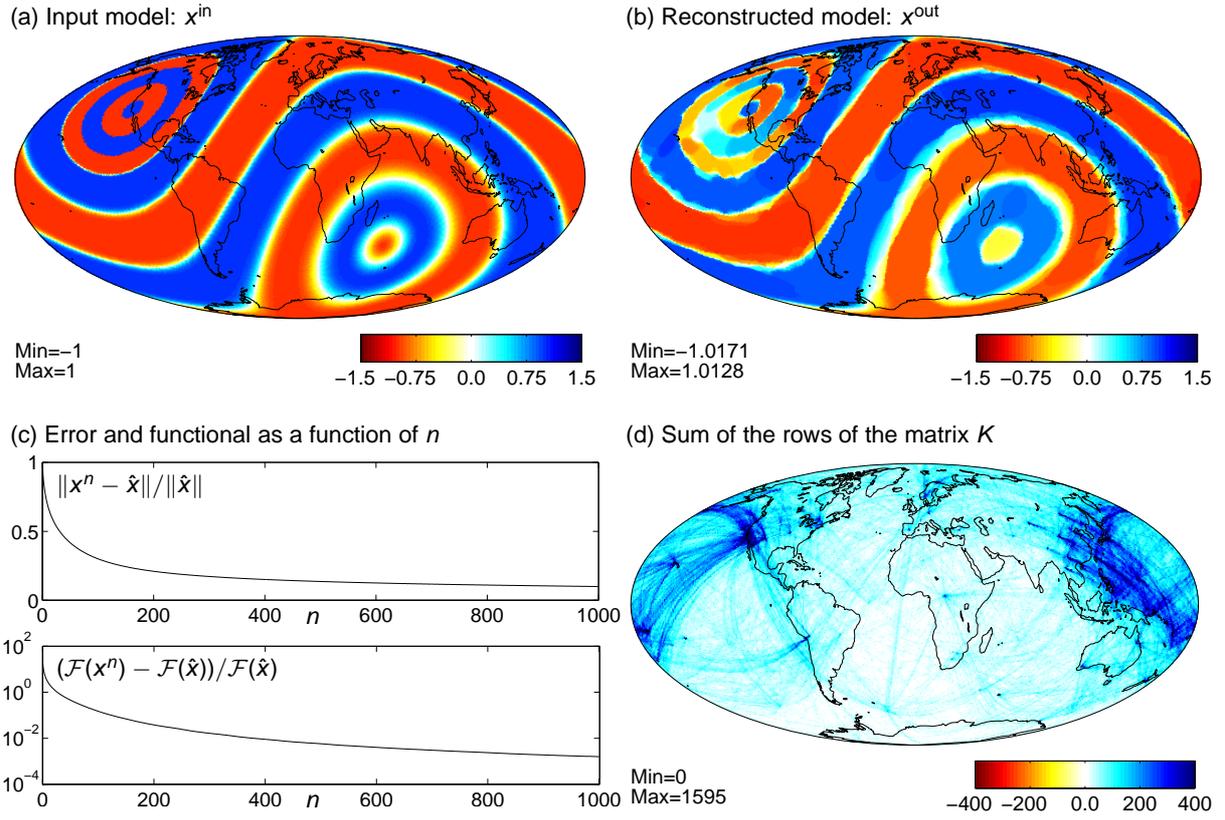}}
\caption{Synthetic seismic tomography experiment using total variation penalty (see section \ref{examplesection}): (a) input model with both sharp and smooth edges between zones of constant model value; (b) reconstruction from $8490$ noisy data with $1000$ iterations of algorithm (\ref{scaledalg}) and $A=\mathrm{grad}$;  (c) evolution of the distance to the limit model and of the functional to the limit value (here $\hat x$ is obtained from $100000$ iterations of the same algorithm; $\hat x$ is not equal to $x^\mathrm{in}$); (d) Sum of the rows of the matrix $K$ to indicate the illumination of the sphere by the rays in the data set.}\label{seismicfigure}
\end{figure}

Algorithm (\ref{scaledalg}) is applied to a stylized problem in seismic tomography.
We try to reconstruct a simple synthetic $2D$ input model $x^\mathrm{in}$ defined on the sphere
(see figure \ref{seismicfigure}.a) from $8490$ data. The input model has a number of zones of constant value with either a sharp edge or a smooth edge in between.
The model space has dimension $98304$. The data $y$ are found from $8490$ seismic surface rays that criss-cross the globe (these correspond to actual earthquakes and seismic stations \cite{Trampert2001}) and that
make up the rows of a matrix $K$ (see figure \ref{seismicfigure}.d). More precisely,
synthetic data $y$ are constructed through the formula
$y=Kx^\mathrm{in}+\epsilon$, where $\epsilon$ is gaussian noise of magnitude
$\|\epsilon\|=0.1\times\|Kx^\mathrm{in}\|$ (i.e. 10\% noise). The aim is now
to reconstruct $x^\mathrm{in}$ as well as possible from $y$ by imposing a total variation penalty in cost function (\ref{functional}). In other words we will look for the minimizer of function (\ref{functional}) where $K$ and $y$ are given and where we choose $A=\mathrm{grad}$.

Algorithm (\ref{scaledalg}), with $\tau=0.99/\|K\|^2$ and $\sigma=0.99/\|A\|^2$ and $1000$ iterations was used to produce a reconstruction: $x^\mathrm{out}=x^{1000}$ (in about $20$ seconds of computer time). The penalty parameter $\lambda$ was chosen to fit the data to the level of the noise: $\|Kx^\mathrm{out}-y\|=\|\epsilon\|$.

The original model and its reconstruction are shown in figure \ref{seismicfigure}, panels (a) and (b). The total variation penalty results in a piece-wise constant output model. Sharp edges (e.g. near North America) are reasonably well resolved given the small amount of data available. In panel (c) the distance of iterate $x^n$ to a reference minimizer $\hat x$ is shown. The residual error (with respect to $\hat x$) after $1000$ iterations is 10\%, and the functional attains about $3$ correct decimal compared to the `true' minimal value.

\section{Acknowledgements}

I.L. is a research associate of the F.R.S.-FNRS  (Belgium).
Part of this research was done while the authors were at CAMP
group of the Vrije Universiteit Brussel and was supported by
VUB GOA-062 and by the FWO-Vlaanderen grant G.0564.09N. The
authors would like to thank Antonin Chambolle for sending them
\cite{Chambolle.Pock2010} and for constructive comments, Frederik Simons for providing them with the list of source-receiver positions used in the synthetic example of section \ref{examplesection} and the two anonymous referees for their valuable remarks.


\end{document}